\documentclass[12pt]{article}
\usepackage[T1]{fontenc}
\usepackage[utf8]{inputenc}
\usepackage{lmodern}
\usepackage{amsmath,amssymb,amsthm}
\usepackage{graphicx}
\usepackage{hyperref}
\usepackage{bm}
\usepackage{microtype}
\usepackage[margin=1in]{geometry}

\title{Integrable nonlinear oscillators with polynomial invariants:\\
construction, Poincar\'e geometry, and an analytic stability boundary}
\author{Johannes Hagel\\[4pt]
Alexander von Humboldt Gymnasium, Neuss\\
Bergheimer Stra\ss e 233, 41464 Neuss, Germany\\
\texttt{johannes.hagel@gmail.com}}
\date{}

\begin{document}
\maketitle

\begin{abstract}
Starting from the nonlinear ODE $z'' + f(t)\,z + g(t)\, z^{m}=0$ with $m>1$, we show that after a suitable normal-form reduction of any Hill equation one may, without loss of generality, fix the linear part as $f(t)\equiv \omega^{2}$ (with $\omega>0$ constant). For the class $z''+\omega^{2}z+g(t)\, z^{m}=0$ with $m>1$, our goal is to compile a catalogue of all possible integrable cases. We restrict attention to integrals that are polynomial in the variables $z$ and $p=z'$. The Hamiltonian does not provide such an integral because it is explicitly time dependent. Instead, we search for invariants that are quadratic in $p=z'$. We show that such invariants exist precisely when $\alpha_2(t):=g(t)^{-2/(m+3)}$ satisfies the linear third-order ODE $\alpha_2''' + 4\omega^2 \alpha_2'=0$. This yields the three-parameter solution $g(t)=[a_0+a_1\cos(2\omega t)+a_2\sin(2\omega t)]^{-(m+3)/2}$. For $m=2$ this reproduces the trigonometric structure with exponent $-5/2$ found in Hagel--Bouquet (1992). In addition we present a detailed stability analysis based on the invariant using Poincar\'e sections and find full agreement with numerical simulations.
\end{abstract}

\section{Introduction}
The aim of this paper is to assemble a systematic catalogue of all integrable cases of the
nonlinear equation

\begin{equation*}
z'' + f(t) z + g(t) z^m = 0, \quad m > 1.
\end{equation*}

We seek first integrals in polynomial form in the variables $z$, $p = z'$, and $t$. The associated
Hamiltonian is not suitable as a conserved integral, because it is explicitly time dependent.
Instead, we search specifically for invariants that are quadratic in $p = z'$, which then permit
a classification of integrable cases.

The nonlinear equation above can be analysed, whenever $f(t)$ is periodic with period $T$, by applying the Courant--Snyder transformation \cite{courant1958} to the \emph{linear part} $z''+f(t)z=0$. In particular,

\begin{equation*}
y = \frac{z}{\sqrt{\beta(t)}}, \qquad \Phi(t) = \int_0^t \frac{d\tau}{\beta(\tau)}, \qquad \omega = \frac{\Phi(T)}{2\pi},
\end{equation*}

(where $w(t)=\sqrt{\beta(t)}$ denotes the periodic envelope of the Hill equation $z'' + f(t)z=0$ as
long as $z(t)$ remains bounded) maps the linear part to the harmonic oscillator

\begin{equation*}
\frac{d^2 y}{d \Phi^2} + \omega^2 y = 0.
\end{equation*}

The nonlinear contribution preserves its form under this change of variables, so that the full equation takes the normal form

\begin{equation*}
z'' + \omega^2 z + g(t) z^m = 0.
\end{equation*}

The structure of the paper is as follows. Section 2 defines the dynamical system in the normal form, introduces the quadratic ansatz for an invariant, and derives the conditions on $g(t)$ under which integrability holds. Section 3 analyzes the invariant curve and establishes an exact analytical stability boundary, which is found to be in complete agreement with numerical simulations. Section 4 summarizes the results and outlines perspectives for further generalizations.
\section{General formulation, integrability ansatz and example}

\subsection{Problem class and objectives}

We first consider the class of explicitly time-dependent nonlinear second-order differential equations
\begin{equation}
z'' + f(t) z + g(t) z^m = 0, \quad m \in \mathbb{N}, \, m \geq 2. \tag{2.1}
\end{equation}
Although this equation can be derived from a Hamiltonian function
\begin{equation}
H(z,p,t) = \frac{p^2}{2} + \frac{f(t)}{2} z^2 + \frac{g(t)}{m+1} z^{m+1}, \tag{2.2}
\end{equation}
it does not provide a conserved integral since
\(\partial H / \partial t \neq 0\). Here $f,g:\mathbb{R}\to \mathbb{R}$ are continuous and, without loss of generality, assumed $2\pi$-periodic. 
Our aim is to identify within class (2.1) special integrable systems. 
We search for first integrals in quadratic polynomial form in $p = z'$:
\begin{equation}
I(z,p,t) = a_0(z,t) + a_1(z,t) p + a_2(z,t) p^2, \qquad \frac{dI}{dt} = 0. \tag{2.3}
\end{equation}

As shown in Chapter 3 of \cite{hagel2015}, condition (2.3) leads to a system of PDEs for the coefficients $a_0,a_1,a_2$ of the form
\begin{align}
- a_1(z,t)\,[z f(t) + z^m g(t)] + \frac{\partial a_0}{\partial t} &= 0, \tag{2.4}\\
-2 a_2(z,t)\,[z f(t) + z^m g(t)] + \frac{\partial a_0}{\partial z} + \frac{\partial a_1}{\partial t} &= 0, \tag{2.5}\\
\frac{\partial a_1}{\partial z} + \frac{\partial a_2}{\partial t} &= 0, \tag{2.6}\\
\frac{\partial a_2}{\partial z} &= 0. \tag{2.7}
\end{align}

The general solutions of (2.5)–(2.7) yield for $a_0,a_1,a_2$ the following expressions:
\begin{align}
a_2 &= \alpha_2(t), \tag{2.8}\\
a_1 &= - \alpha_2'(t) z + \alpha_1(t), \tag{2.9}\\
a_0 &= z^2 \alpha_2(t) f(t) + \frac{2}{m+1} z^{m+1} \alpha_2(t) g(t) - z \alpha_1'(t) + \tfrac{1}{2} z^2 \alpha_2''(t) + \alpha_0(t). \tag{2.10}
\end{align}

Inserting these into (2.4) and comparing coefficients of $z$ yields the following cases:

\paragraph{(i) $m=2$:}
\begin{align}
\alpha_1'' + f(t) \alpha_1 &= 0, \tag{2.11}\\
\alpha_2''' + 4 f(t)\alpha_2' + 2 \alpha_2 f'(t) - 2 \alpha_1 g(t) &= 0, \tag{2.12}\\
5 g(t)\alpha_2' + 2 \alpha_2 g'(t) &= 0. \tag{2.13}
\end{align}

\paragraph{(ii) $m>2$:}
\begin{align}
\alpha_1'' + f(t) \alpha_1 &= 0, \tag{2.14}\\
\alpha_2''' + 4 f(t)\alpha_2' + 2 \alpha_2 f'(t) - 2 \alpha_1 g(t) &= 0, \tag{2.15}\\
-\alpha_1 g(t) &= 0, \tag{2.16}\\
(m+3) g(t)\alpha_2' + 2 \alpha_2 g'(t) &= 0. \tag{2.17}
\end{align}

\subsection{The case $f(t) = \omega^2$}

From Chapter 1 it was established that in (2.1) we may set $f(t)=\omega^2$ without loss of generality. 
Then (2.11)–(2.17) simplify to

\paragraph{(i) $m=2$:}
\begin{align}
\alpha_1'' + \omega^2 \alpha_1 &= 0, \tag{2.18}\\
\alpha_2''' + 4\omega^2 \alpha_2' - 2 \alpha_1 g(t) &= 0, \tag{2.19}\\
5 g(t)\alpha_2' + 2 \alpha_2 g'(t) &= 0. \tag{2.20}
\end{align}

\paragraph{(ii) $m>2$:}
\begin{align}
\alpha_1'' + \omega^2 \alpha_1 &= 0, \tag{2.21}\\
\alpha_2''' + 4\omega^2 \alpha_2' - 2 \alpha_1 g(t) &= 0, \tag{2.22}\\
-\alpha_1 g(t) &= 0, \tag{2.23}\\
(m+3) g(t)\alpha_2' + 2 \alpha_2 g'(t) &= 0. \tag{2.24}
\end{align}

For $m=2$, (2.18) yields
\begin{equation}
\alpha_1(t) = \tfrac{1}{2} C_1 \cos(\omega t) + \tfrac{1}{2} C_2 \sin(\omega t), \quad C_{1,2}\in \mathbb{R}. \tag{2.25}
\end{equation}
Eq. (2.20) implies
\begin{equation}
g(t) = \alpha_2(t)^{-5/2}. \tag{2.26}
\end{equation}
Thus $\alpha_2(t)$ satisfies
\begin{equation}
\alpha_2''' + 4 \omega^2 \alpha_2' - [C_1 \cos(\omega t) + C_2 \sin(\omega t)] \alpha_2^{-5/2} = 0. \tag{2.27}
\end{equation}
Hence, in principle, there exists a five-parameter family of $g(t)$ leading to integrable systems of the form
\begin{equation}
z'' + \omega^2 z + g(t) z^2 = 0. \tag{2.28}
\end{equation}

Since no closed form of (2.27) is known, we restrict to the case $C_1=C_2=0$, which renders (2.27) linear with constant coefficients. Its general solution is three-parametric:
\begin{equation}
\alpha_2(t) = A + B \cos(2\omega t) + C \sin(2\omega t). \tag{2.29}
\end{equation}
Thus
\begin{equation}
g(t) = [A + B \cos(2\omega t) + C \sin(2\omega t)]^{-5/2}. \tag{2.30}
\end{equation}

For $m>2$, equation (2.23) implies $\alpha_1=0$ (since $g(t)\neq 0$). Eq. (2.22) then yields
\begin{equation}
\alpha_2(t) = A + B \cos(2\omega t) + C \sin(2\omega t). \tag{2.31}
\end{equation}
From (2.24) we find after integration
\begin{equation}
\alpha_2(t) = g(t)^{-\tfrac{2}{m+3}}. \tag{2.32}
\end{equation}
Hence the integrable cases are given by
\begin{equation}
g(t) = [A + B \cos(2\omega t) + C \sin(2\omega t)]^{-(m+3)/2}, \tag{2.33}
\end{equation}
leading to integrable systems of the form
\begin{equation}
z'' + \omega^2 z + g(t) z^m = 0. \tag{2.34}
\end{equation}

\subsection{Integral for $m=2$ and demonstration example}

For $m=2$, a three-parameter form of $\alpha_2(t)$ is given by
\begin{equation}
\alpha_2(t) = A + B \cos(2\omega t) + C \sin(2\omega t), \qquad g(t)=\alpha_2(t)^{-5/2}. \tag{2.35}
\end{equation}
Eq. (2.34) then admits the quadratic invariant
\begin{align}
I(z,p,t) &= \alpha_2(t) p^2 - \alpha_2'(t) z p \notag \\
&\quad + \Big( \omega^2 \alpha_2(t) + \tfrac{1}{2}\alpha_2''(t) \Big) z^2 + \tfrac{2}{3}\alpha_2(t) g(t) z^3. \tag{2.36}
\end{align}
A Poincar\'e section at times $t_k = k\pi/\omega$ then yields an algebraic curve $p=p(z)$.

\begin{figure}[t]
 \centering
 \includegraphics[width=\linewidth]{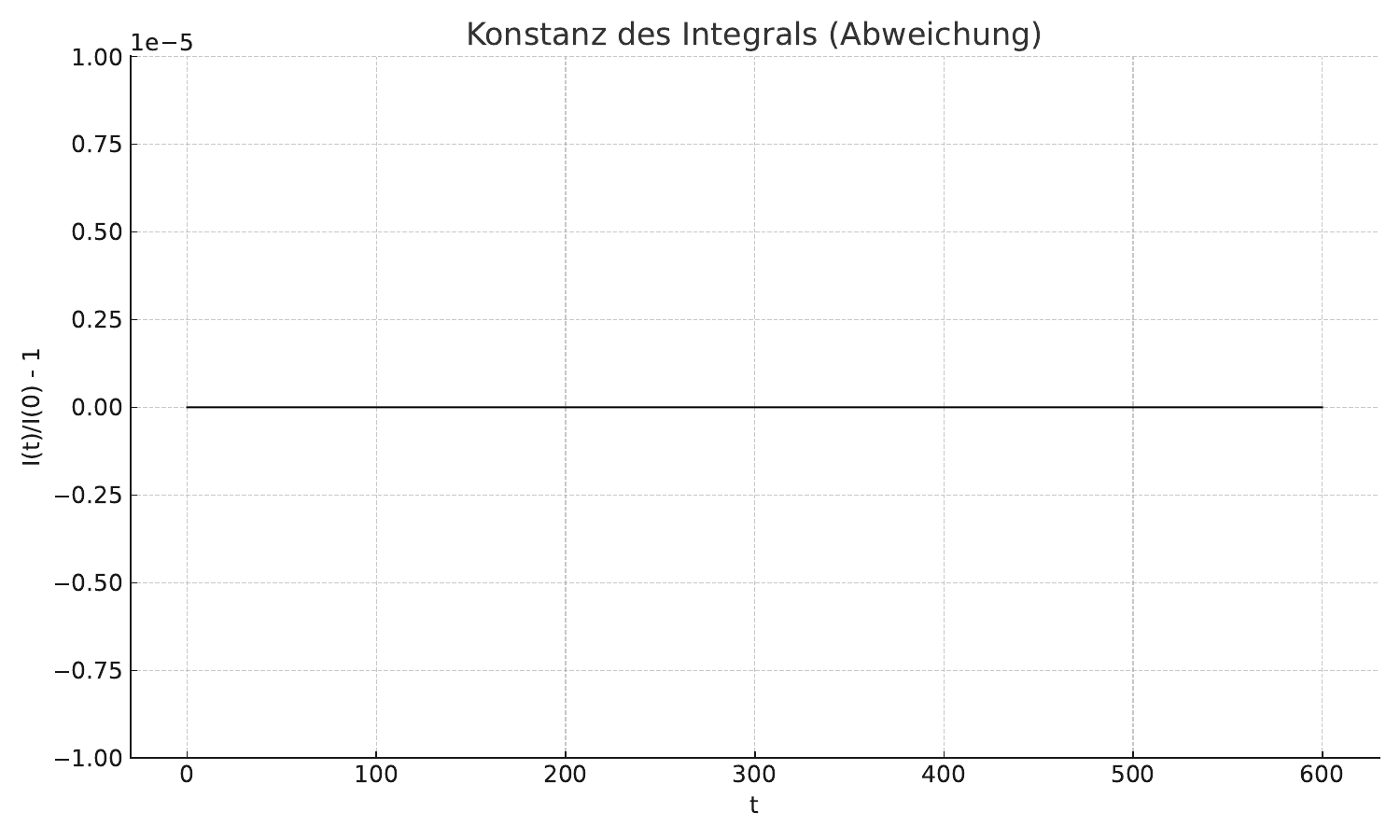}
 \caption{Constancy of the invariant: plot of $I(t)/I(0)-1$ over $0<t<600$ (parameters as in text). The range shown is $[-10^{-5},+10^{-5}]$.}
 \label{fig:constancy}
\end{figure}

As an example, choose $A=1.3$, $B=0.9$, $C=0$, and $\omega=1$. 
We integrate (2.34) with a Runge--Kutta method (step size $h$), initial conditions $z(0)=0.1$, $z'(0)=0$, over $0<t<600$. 
We then plot $I(t)/I(0)-1$. Figure~\ref{fig:constancy} shows the result in the range $[-10^{-5},+10^{-5}]$: the curve runs numerically flat at $0$, as expected, confirming integrability of (2.34) under assumption (2.35).

\subsection{Phase curve at the Poincar\'e section $t_k = k\pi/\omega$}

For $C=0$, at $t_k = k\pi/\omega$ we have $\alpha_2(t_k)=A+B$, $\alpha_2'(t_k)=0$, $\alpha_2''(t_k)=-4\omega^2 B$. With $I_0=I(z(0),z'(0),0)$ the phase curve becomes
\begin{equation}
(A+B) p^2 + \omega^2 (A-B) z^2 + \tfrac{2}{3} (A+B)^{-3/2} z^3 = I_0. \tag{2.37}
\end{equation}
This is quadratic in $p$, yielding
\begin{equation}
p(z) = \pm \sqrt{\frac{I_0 - \omega^2(A-B)z^2 - \tfrac{2}{3}(A+B)^{-3/2} z^3}{A+B}}. \tag{2.38}
\end{equation}

\begin{figure}[t]
 \centering
 \includegraphics[width=0.9\linewidth]{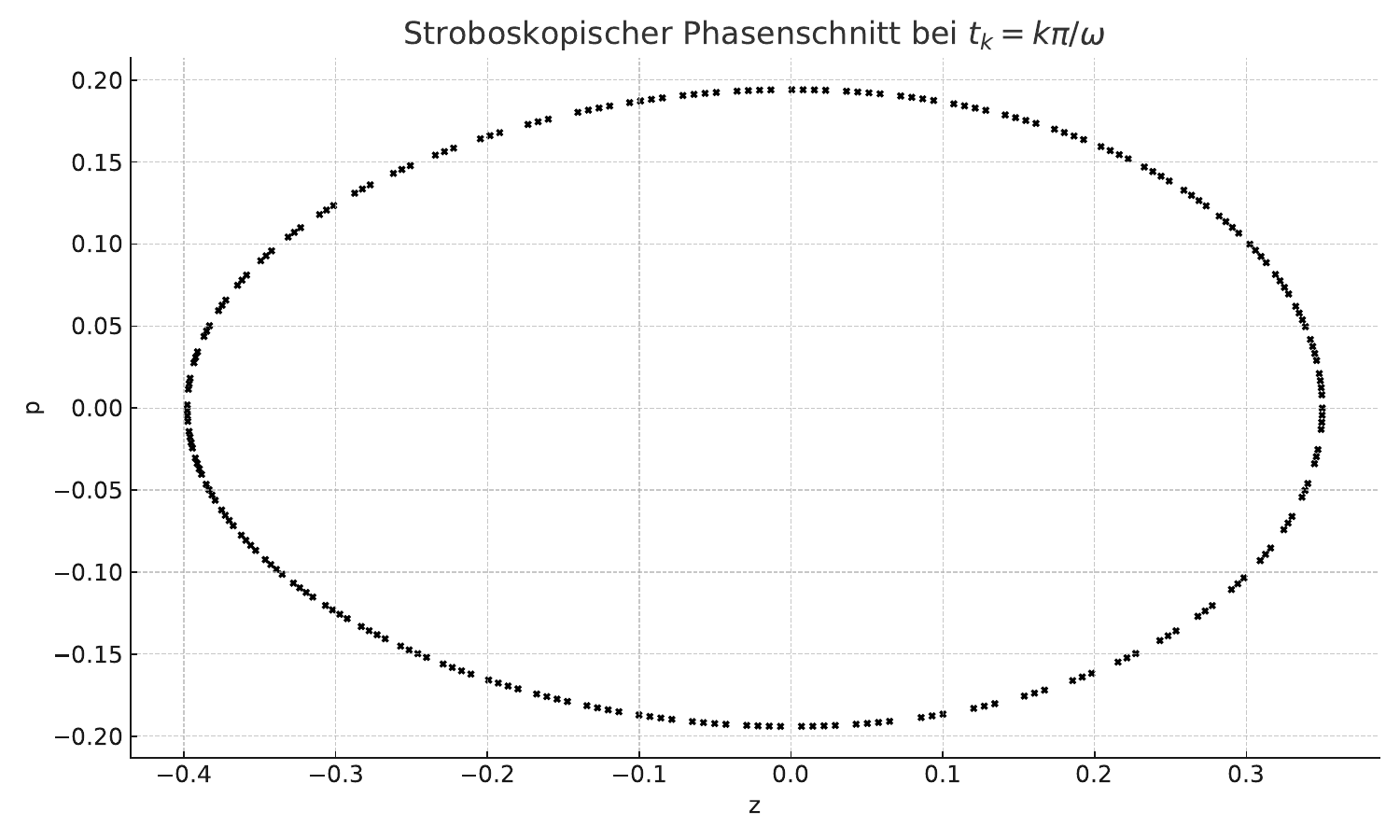}
 \caption{Poincar\'e phase section at $t_k = k\pi/\omega$: stroboscopic points $(z(t_k),p(t_k))$ together with the analytic curve given by (2.38).}
 \label{fig:poincare}
\end{figure}

\section{Stability analysis by means of the invariant curve}

We now investigate the stability boundary of equation (2.34) with specialization (2.35). For this purpose we consider the stroboscopic section at times $t_k=k\pi/\omega$, since at these points $\dot{\alpha}_2(t_k)=0$ and the invariant assumes a particularly simple form.

Set $\alpha=A+R\cos(2\omega t-\varphi)$ with $R=\sqrt{B^2+C^2}$ and suitable phase $\varphi$. The associated stroboscopic curve is then
\begin{equation}
(A+R) p^2 + \omega^2 (A-R) z^2 + \tfrac{2}{3}(A+R)^{-3/2} z^3 = I_0, \tag{3.1}
\end{equation}
where $I_0$ is determined by the initial condition $z(0)=z_0,\,p(0)=0$:
\begin{equation}
I_0 = \omega^2 (A-R) z_0^2 + \tfrac{2}{3}(A+R)^{-3/2} z_0^3. \tag{3.2}
\end{equation}

The shape of the admissible set $p^2\ge 0$ changes precisely when the discriminant of the cubic polynomial in $z$ vanishes. This yields the critical integral value
\begin{equation}
I_{0,\mathrm{crit}} = \tfrac{4}{27}\,[\omega^2(A-R)]^3 \Big(\tfrac{2}{3}(A+R)^{-3/2}\Big)^2
= \tfrac{1}{3}\,\omega^6 (A^2-R^2)^3. \tag{3.3}
\end{equation}
Setting $I_0=I_{0,\mathrm{crit}}$ in the relation for $I_0(z_0)$ then gives the critical initial amplitude
\begin{equation}
z_{\mathrm{crit}}(\omega) = \frac{\omega^2}{2}\,(A-R)(A+R)^{3/2}. \tag{3.4}
\end{equation}

\medskip\noindent\textbf{Theorem 1.} 
For initial conditions $p(0)=0$, $z(0)=z_0$ the motion is bounded as long as $z_0<z_{\mathrm{crit}}(\omega)$. Above this threshold the stroboscopic curve opens and the trajectory leaves the bounded domain.

\medskip
For the parameters of Section~2.3 ($A=1.3$, $B=0.9$, $C=0$) we obtain $R=0.9$ and therefore
\begin{equation}
z_{\mathrm{crit}}(\omega) \approx 0.5\,\omega^2\cdot 0.4\cdot (2.2)^{3/2}. \tag{3.5}
\end{equation}
For $\omega=1.23$ we find $z_{\mathrm{crit}}\approx 0.99$, so that the chosen initial condition $z(0)=0.35$ of Section~2.3 lies well inside the stable regime.

\begin{figure}[t]
 \centering
 \includegraphics[width=0.9\linewidth]{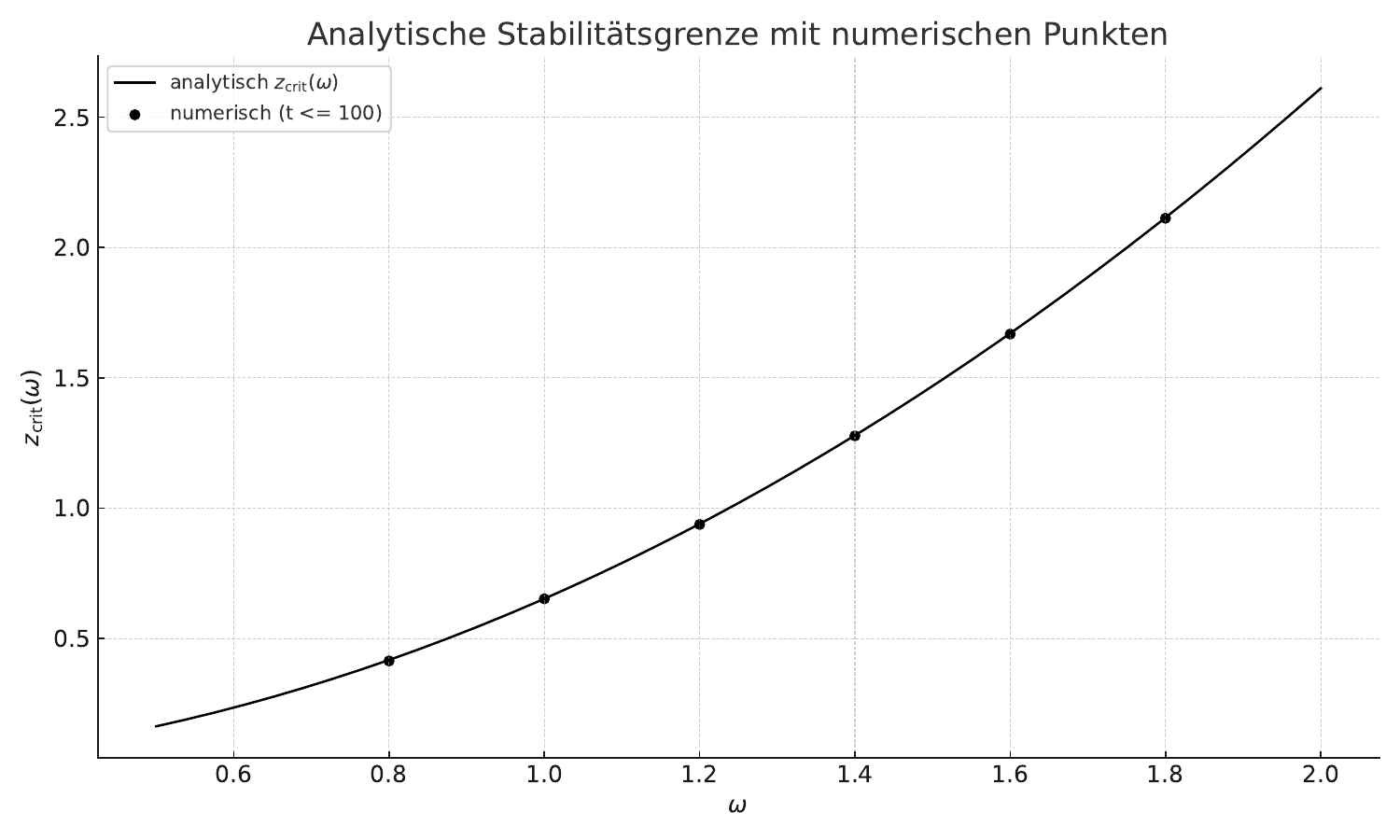}
 \caption{Analytical stability boundary $z_{\mathrm{crit}}(\omega)$ from (3.4) for the example parameters $A=1.3$, $B=0.9$, $C=0$. Black dots: numerical results.}
 \label{fig:stab}
\end{figure}

To test the theory we performed numerical integrations for $\omega=0.8,1.0,1.2,1.4,1.6,1.8$ with increasing amplitude in steps of $\Delta z_0=0.02$, marking the last initial value $z_0$ that still led to a bounded solution $z(t;z_0)$ with a black dot. Figure~\ref{fig:stab} shows that these results for the stability boundary are in perfect agreement with the analytic curve (3.4)--(3.5).

\begin{figure}[t]
 \centering
 \includegraphics[width=0.9\linewidth]{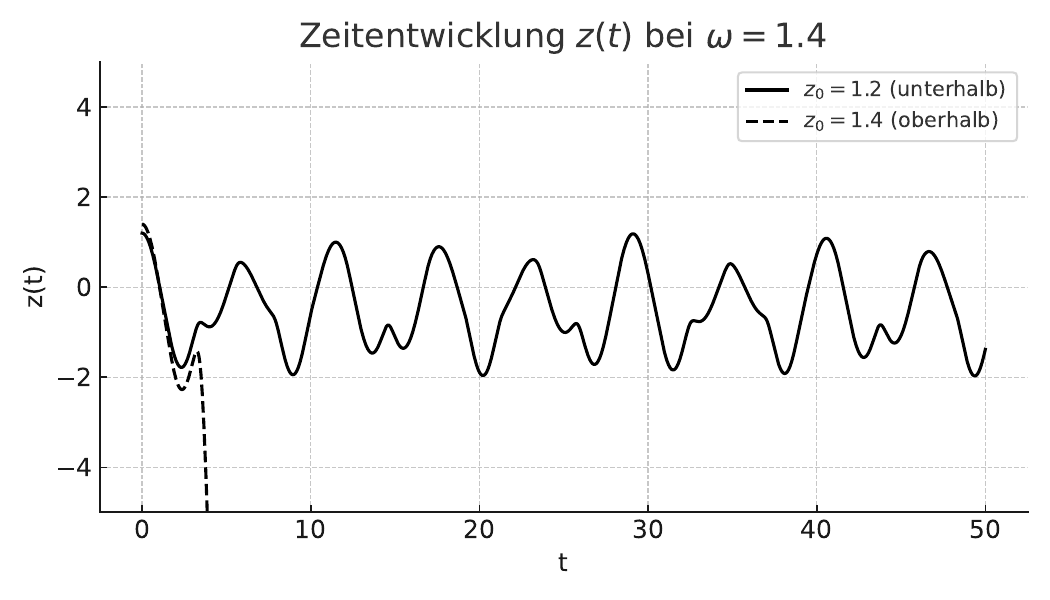}
 \caption{Time evolution $z(t)$ for $\omega=1.4$. Solid line: $z(0)=1.2$ (below the stability threshold, bounded). Dashed line: $z(0)=1.4$ (above the threshold, unbounded). The vertical axis is scaled to $[-5,5]$.}
 \label{fig:time}
\end{figure}

Figure~\ref{fig:time} displays one bounded and one unbounded solution $z(t)$, respectively for $z_0=1.2$ and $z_0=1.4$.

\section{Conclusions}

We developed a systematic approach to identify integrable cases of
\begin{equation*}
z''+\omega^{2}z+g(t)\,z^{m}=0,\qquad m>1,
\end{equation*}
by constructing invariants that are quadratic in $p=z'$ and deriving the corresponding conditions on $g(t)$ under which integrability holds.

For $m=2$ a five-parameter family arises in principle, which reduces in a special case to a three-parameter family with trigonometric structure. For $m>2$ only three-parameter families occur. The conservation of the invariants was verified numerically; deviations were observed only at the level of round-off. Poincar\'e sections illustrated the closed invariant surfaces. Moreover, we obtained an exact analytical stability boundary whose graph is in complete agreement with numerical experiments.

These results provide a coherent picture of the integrable classes of nonlinear oscillators of this form and constitute a starting point for future investigations towards generalized invariants and higher nonlinearities.

\section*{Acknowledgments}
The author thanks in particular S.~Bouquet for earlier joint work that initiated this investigation. Special thanks are due to Gilbert Guignard for intensive joint research during the years 1983--1992. The author is also grateful to colleagues in accelerator physics for valuable discussions on integrability and stability in nonlinear systems. Finally, the author thanks ChatGPT (OpenAI) for assistance with editing and structuring the manuscript. The comments of anonymous referees are also gratefully acknowledged.


\begin{thebibliography}{9}

\bibitem{hagel-bouquet-1992}
J.~Hagel, S.~Bouquet, \emph{Integrals for a special class of second order equations with a quadratic nonlinearity}, CERN SL/92-52 (AP), 1992.

\bibitem{hagel2015}
J.~Hagel, \emph{A new method to construct integrable approximations to nearly integrable systems in Celestial Mechanics: application to the Sitnikov problem}, Celest.\ Mech.\ Dyn.\ Astron.\ \textbf{122} (2015), 101--116.

\bibitem{courant1958}
E.~D.~Courant, H.~S.~Snyder, \emph{Theory of the Alternating-Gradient Synchrotron}, Ann.\ Phys.\ \textbf{3} (1958), 1--48.

\bibitem{guignard1978}
G.~Guignard, \emph{The general theory of all sum and difference resonances in a three-dimensional magnetic field in a synchrotron}, CERN Report 78-11, 1978.

\bibitem{guignard-verdier-1980s}
G.~Guignard, A.~Verdier, \emph{Analytical theory and correction of nonlinear resonances in synchrotrons}, CERN Internal Reports, 1980s.

\bibitem{guignard-hagel-1986}
G.~Guignard, J.~Hagel, \emph{Sextupole Correction and Dynamic Aperture: Numerical and Analytical Tools}, Particle Accelerators \textbf{18} (1986) 129--165. (CERN-LEP-TH-85-3).

\end{thebibliography}
\end{document}